\documentclass{amsart}

\usepackage{amsthm,amsfonts,amsmath,amssymb,latexsym,enumitem}
\usepackage{upref,amssymb,eucal,ae}
\usepackage[all,cmtip]{xy}

\newtheorem{theorem}{Theorem}

\newtheorem{corollary}[theorem]{Corollary}

\newenvironment{remark}{\medskip \refstepcounter{theorem}
\noindent  {\bf Remark \thetheorem}.\rm}{\,}

\makeatletter
\def\Ddots{\mathinner{\mkern1mu\raise\p@
\vbox{\kern7\p@\hbox{.}}\mkern2mu
\raise4\p@\hbox{.}\mkern2mu\raise7\p@\hbox{.}\mkern1mu}}
\makeatother

\newcounter{spthe}

\def\<{\langle}
\def\a{\alpha}
\def\>{\rangle}

\def\tm{\tilde{M}}
\def\tn{\tilde{n}}
\def\tg{\tilde{g}}
\def\mb#1{{\mathbb #1}}
\def\mc#1{{\mathcal #1}}

\def\BOne{{\mathchoice {\rm 1\mskip-4mu l} {\rm 1\mskip-4mu l}
                          {\rm 1\mskip-4.5mu l} {\rm 1\mskip-5mu l}}}

\begin{document}
\title[Minimally embedded Riemann surfaces in $\mb{S}^3$]{Minimally  
embedded Riemann surfaces in $\mb{S}^3$ and the conformal deformation of
their metrics}
\author{Santiago R. Simanca}
\email{srsimanca@gmail.com}

\begin{abstract}
We prove that if $f_g: (\Sigma,g) \rightarrow (\mb{S}^{2+p},\tg)$ is a 
smooth minimal isometric embedding of a Riemannian surface $(\Sigma,g)$, and 
$[0,1]\ni t \rightarrow g_t$ 
is a path of area preserving conformal deformations of $g$ on 
the embedded submanifold $f_g(\Sigma)$, then there exists a path of conformal
diffeomorphism $F_t: (\mb{S}^{2+p}, F_t^*\tg) \rightarrow (\mb{S}^{2+p},\tg)$ 
that starts at $\BOne_{\mb{S}^{2+p}}$, set theoretically fixes $f_g(\Sigma)$ 
for all $t$, and it is such that $F^{*}_t \tilde{g}\mid_{f_g(M)}=g_t$ with 
$f_{g_t}: (\Sigma,g_t) \rightarrow (\mb{S}^{2+p},\tg)$ a path of minimal
embedding deformations of the initial $f_g$. We apply this result to the 
Lawson surface $(\Sigma,g)=(\xi_{k/m,m},g_{\xi_{k/m,m}})$, $m$ any divisor 
of $k>1$, to conclude that if $a:=\mu_{g_{\xi_{k/m,m}}}(\Sigma)$, 
and $[0,1]\ni t \rightarrow g_t$ is a path of area $a$ metrics conformal 
deformations of $g_{k/m,m}$ to a metric $g_a$ of constant scalar curvature 
$s_{g_a}=4\pi \chi(\Sigma)/a$, 
the embedding $f_{g_{\xi_{k/m,m}}}: (\xi_{k/m,m},g_{\xi_{k/m,m}}) \rightarrow 
(\mb{S}^3,\tg)$ has associated minimal isometric conformal deformations 
$f_{g_t}$ to the isometric embedding $f_{g_a}$ of $g_a$,  
in sharp contrast 
with the situation of the standard sphere $(\xi_{0,1},g_{\xi_{0,1}})$ and
Clifford torus $(\xi_{1,1},g_{\xi_{1,1}})$, which
are the only orientable Riemannian surfaces 
of genus $0$ and $1$, respectively, 
isometrically embedded into $(\mb{S}^3,\tg)$ as minimal surfaces. 
If $\sigma^2(\Sigma):=\sup_{[g]\in \mc{C}(\Sigma)}(4\pi 
\chi(\Sigma))^2/\left(\frac{1}{4}\inf_{g\in [g]}\mc{W}_{f_g}(\Sigma)\right)$, 
$\mc{W}_{f_g}(\Sigma)$ the Willmore energy of the isometric 
embedding $f_g:(\Sigma,g)\rightarrow (\mb{S}^{\tn},\tg)$ 
and $\mc{C}(\Sigma)$ the space of all conformal classes of metrics on 
$\Sigma$, if $\Sigma^k$ is orientable then 
$(4\pi \chi(\Sigma))^2/\left(\frac{1}{4}\mc{W}_{f_g}(\Sigma) \right) \leq 
\sigma^2(\Sigma)=(4\pi \chi(\Sigma))^2/\left(\frac{1}{4}\mc{W}_{f_{
g_{\xi_{k,1}}}}(\Sigma) \right)$, and for $k\neq 1$, the equality is achieved 
by $f_g$ if, and only if, $[g]=[g_{\xi_{k,1}}]$ and $f_g$ 
is conformally equivalent to $f_{g_{\xi_{k,1}}}$, while for $k=1$, 
$\chi(\Sigma)=0$ so the equality is always achieved, but then 
$\mc{W}_{f_g}(\Sigma)=\mc{W}_{f_{\xi_{1,1}}}(\Sigma)$ if, and only if,
$[g]=[g_{\xi_{1,1}}]$ and $f_g$ is conformally equivalent to $f_{\xi_{1,1}}$. 
\end{abstract}
\maketitle

\section{Isometric embeddings into the standard sphere} \label{s1}
If $(M^{n\geq 2},g)$ is a closed Riemannian manifold, and 
\begin{equation} \label{emb}
f_g: (M,g) \rightarrow (\tm^{\tn}, \tg)
\end{equation}
is an isometric embedding of it into some Riemannian background 
$(\tm^{\tn}, \tg)$, the scalar curvature $s_g$ of $g$ relates to 
extrinsic quantities of $f_g$ by 
\begin{equation} \label{sce}
s_g  = \sum K^{\tg}(e_i,e_j) +\tg(H,H)- \tg(\a,\a) \, .
\end{equation}
Here, $\{ e_i\}$ is an orthonormal tangent frame of $f_g(M)$, $K^{\tg}$ is the 
sectional curvature of the metric $\tg$, and $\tg(H,H)$ and 
$\tg(\alpha,\alpha)$ are the squared $\tg$-norms of the mean curvature 
vector $H:=H_{f_g}$, and second 
fundamental form $\alpha:= \alpha_{f_g}$ tensors of $f_g$, 
respectively. By integration of these three extrinsic functions 
with respect to the intrinsic volume measure of $g$, we obtain 
extrinsic functionals $\Theta_{f_g}(M)$, $\Psi_{f_g}(M)$ and $\Pi_{f_g}(M)$, 
respectively, and we have that 
\begin{equation} \label{eq11} 
\begin{array}{rcl} 
{\displaystyle \mc{S}_{g}(M):= \int_{f_g(M)} s_{g} d\mu_{g} } 
& =  &  
 \Theta_{f_g}(M) + \Psi_{f_g}(M)  
-{\displaystyle \Pi_{f_g}(M)} 
\\ & = & {\displaystyle  \mc{W}_{f_g}(M) - \mc{D}_{f_g}(M)} 
 \, ,
\end{array}
\end{equation}
where 
\begin{equation} \label{eq12}
\begin{array}{rcl}
\mc{W}_{f_g}(M) & := & {\displaystyle
\frac{n}{n-1}\Theta_{f_g}(M)+\Psi_{f_g}(M) }\, ,\vspace{1mm}\\
\mc{D}_{f_g}(M) & := & {\displaystyle
\frac{1}{n-1}\Theta_{f_g}(M)+\Pi_{f_g}(M)}  \, .
\end{array}
\end{equation}
The particular combinations of 
extrinsic functionals makes of 
$\mc{W}_{f_g}(M)$ and $\mc{D}_{f_g}(M)$ 
intrinsically defined functionals on the 
space of metrics on the manifold in a fixed conformal class,  
a convenient fact for the geometric analysis of $M$ if the background 
$(\tm,\tg)$ were to carry isometric embeddings of all the metrics of $M$. 

By the Nash isometric embedding theorem \cite{nash}, any Riemannian manifold 
$(M^n,g)$ can be isometrically embedded into a standard sphere 
$(\mb{S}^{\tn},\tg) \hookrightarrow (\mb{R}^{\tn+1},\| \, \, \|^2)$ of large
but fixed dimension $\tn=\tn(n)$. We thus 
identify the open cone $\mc{M}(M)$ of smooth Riemannian metrics on 
$M$ with the space of their isometric embeddings into this sphere background,  
$$ 
\mc{M}_{\mb{S}^{\tn}}(M)=\{f_g(M)\subset \mb{S}^{\tn}:\; f_g: (M,g) 
\rightarrow (\mb{S}^{\tn},\tg) \; \text{isometric embedding}\}\, ,   
$$ 
and the deformations of metrics in $\mc{M}(M)$ with the 
isotopic Palais' deformations of their isometric embeddings \cite{pal}.
If $f_g\in \mc{M}_{\mb{S}^{\tn}}(M)$, 
the exterior scalar curvature on $f_g(M)$ is 
$\sum K^{\tg}(e_i,e_j)=n(n-1)$, and since there exist volume preserving 
deformations of $g$ in the conformal class $[g]$, 
the isometric embedding $f_g$ of $g$ must be such that 
$\| H_{f_g}\|^2$ is a constant function on $f_{g}(M)$ \cite[Theorem 6]{scs}. 

\subsection{Principal example: Lawson's minimal surfaces immersed
into $\mb{S}^3$} 
If $\Sigma=\Sigma^k$ is a a surface of topological genus $k$, 
an immersion $f : \Sigma \rightarrow \mb{S}^m$ into the standard sphere 
fixes a conformal class on $\Sigma$, who inherits an intrinsic metric 
representative induced by restriction of the metric of the background sphere 
to $f(\Sigma)$. If 
$\mc{I}_{\mb{S}^m}(\Sigma^k)=\{f: \; f : \Sigma^k \rightarrow \mb{S}^m 
\text{\mbox{} is an immersion}\}$ is the set of all possible immersions of 
$\Sigma$ into $\mb{S}^m$, other than for 
$\Sigma =\mb{P}^2(\mb{R})$, $\mc{I}_{\mb{S}^3}(\Sigma^k)$ is nonempty 
and has distinguished elements in it given by surfaces 
$\xi_{m',m}$, $\tau_{m',m}$, and $\eta_{m',m}$ associated to pairs of 
nonnegative integers $m', m$ that determine the surface's topological genus, 
all minimally embedded or immersed into $\mb{S}^3$, discovered 
by Lawson back in the late 1960s \cite[Corollary 1.6, Theorems 2,3,4]{la2}. 
The surfaces $\xi_{m',m}$ are orientable while the surfaces 
$\eta_{m',m}$ are not. Some but not all of the $\tau_{m',m}$s are
nonorientable. The immersions of the $\xi_{m',m}$s are by embeddings.

We denote by $g_{\xi_{m',m}}$ the intrinsic Riemannian metric on the
Riemann surface $\xi_{m',m} \hookrightarrow \mb{S}^3$.
If $m'm=0$, say $m'=0$ and $m>0$, $\xi_{0,m}$ is the standard two sphere in 
$\mb{S}^3$, in which case we associate it to the single pair $(m',m)=(0,1)$;
otherwise, $\xi_{m',m}$ is a Riemann surface of genus $k=m'm$ 
and the embedding into $\mb{S}^3$ is minimal.    
When $k>1$ is not a prime and $1,k\neq m\mid k$, the Riemann surfaces 
$\xi_{k,1}$ and $\xi_{k/m,m}$, both of genus $k$, 
have intrinsic metrics in different conformal classes
distinguished from each other by the symmetries of their fundamental domains. 
We shall single out the surfaces $\xi_{k,1}$, which for $k\geq 1$ 
have equilateral fundamental domain \cite[Proposition 6.1]{la2}, 
and refer to them all (including the sphere $\xi_{0,1}$) as the equilateral 
Lawson Riemann surface of genus $k$. 

The sphere $\xi_{0,1}=\mb{S}^2$ is geodesically embedded
into $\mb{S}^3$ with area $\mu_{g_{\xi_{0,1}}}(\xi_{0,1})=4\pi$.  
The Clifford torus $\xi_{1,1}=\mb{S}^2(1/\sqrt{2})\times \mb{S}^2(1/\sqrt{2})$
is linearly embedded into $\mb{S}^3$ with area 
$\mu_{g_{\xi_{1,1}}}(\xi_{1,1}) = 2\pi^2$.
No 
such an explicit expression for $\mu_{g_{\xi_{k,1}}}(\xi_{k,1})$ exists
when $k>1$, but it is known that 
none of the then metrics $g_{\xi_{k,1}}$ are of constant scalar curvature 
\cite[Proposition 1.5]{la2}, that  $\mu_{g_{\xi_{k,1}}}(\xi_{k,1})< 8\pi$, and 
that $\mu_{g_{\xi_{k,1}}}(\xi_{k,1})\nearrow 8\pi$ as $k \nearrow \infty$ 
\cite{ku,kls}. 

If $\Sigma$ is not orientable, if $k>2$, the minimal surface 
$\eta_{k-1,1}=f_{\eta_{k-1,1}}(\Sigma) \hookrightarrow \mb{S}^3$ has 
topological genus $k$, but in general has self-intersections 
\cite[\S 8]{la2}. 
A nonorientable surface of topological genus $k=2$ is a Klein bottle, and 
Lawson realizes the most symmetric one of them as the associated bipolar 
surface $\Sigma=\tilde{\tau}_{3,1}$ of his manifold $\tau_{3,1}$, which is in 
the conformal class of the square torus $\mb{R}^2/\Gamma$ of lattice 
$\Gamma$ generated by $(\pi,\pi)$ and $(\pi,-\pi)$, and can be explicitly 
described by a doubly periodic immersion into $\mb{S}^3$ 
\cite[\S 7 \& \S 11]{la2}.
The parametrization of $\Sigma=\tilde{\tau}_{3,1}$ can be given then
explicitly \cite[\S 3]{lap}, and 
we obtain a minimal isometric embedding
$f_g: (\tilde{\tau}_{3,1},g_{\tilde{\tau}_{3,1}}) \rightarrow (\mb{S}^4,\tg)
\hookrightarrow (\mb{S}^5, \tg)$, where 
by an explicit calculation, the metric is  
of area $\mu_{g_{\tilde{\tau}_{3,1}}}(\tilde{\tau}_{3,1})= 
6\pi E(2\sqrt{2}/3)$, $E$ the complete elliptic integral.
The surface $\mb{P}^2(\mb{R})$, of genus $1$, can be 
minimally embedded into $\mb{S}^4$ by a map from
the standard sphere $\mb{S}^2(r)$ of radius $r=\sqrt{3/2}$ that is a
$2$-to-$1$ cover of the embedding, with the conformal structures
of base and cover consistent with each other \cite{cdck,si3}, and area
$\mu_{g_{\mb{P}^2}}(\mb{P}^2(\mb{R}))=6\pi$. 

\section{Conformal deformations of isometric embeddedings of surfaces}
\label{s2}
Given the embedding (\ref{emb}) of $(M,g)$ into $(\mb{S}^{\tn},\tg)$,
 suppose that 
\begin{equation} \label{emc}
[0,1]\ni t \rightarrow f_{g_t}: (M,g_t) \mapsto (\mb{S}^{\tn}, \tg)
\end{equation}
is a path of conformally related isometric embedding deformations of 
$f_g=f_{g_t\mid_{t=0}}$ of conformally related Riemannian metrics 
\begin{equation} \label{cd}
[0,1]\ni t \rightarrow g_t =e^{2\psi (t)}g
\end{equation}
on $M$, $\psi(t)$ a path of functions, $\psi(0)=0$.
By the Palais isotopic extension theorem, there exists a smooth one 
parameter family of diffeomorphisms
$$
F_t: \mb{S}^{\tn} \rightarrow \mb{S}^{\tn}
$$
such that $F_t(f_g(x))=f_{g_t}(x)$, and, by restriction, we obtain a 
diffeomorphism 
$$
F_t: (f_g(M),\tg) \rightarrow (f_{g_t}(M), \tg) \, . 
$$
The tensor $F_t^* \tg$ is just the metric $\tg$ on $\mb{S}^{\tn}$ acted on by 
the diffeomorphism $F_t$, and since the $g_t$s are conformal deformations of
$g$, we have that
$$
F_t^* \tg\mid_{f_{g_t}(M)} = e^{2 u(t)(f_g(\, \cdot \,))} \tg \mid_{f_g(M)} =
e^{2 u(t)(f_g(\, \cdot \,))} g  \, ,
$$
where the conformal factor $e^{2 u(t)}$ and that in (\ref{cd}) are related 
to each other by $e^{2\psi(t)(\, \cdot \, )}=e^{2u(t)\circ f_g(\, \cdot \, )}$. 
We extend $u(t)$ conveniently to a function on the whole of
$\mb{S}^{\tn}$, and view the family (\ref{emc}) as the family of conformal 
isometric embeddings 
$$ 
f_{g_t}: (M,e^{2u(t)\circ f_g}g) \hookrightarrow (\mb{S}^{\tn},e^{2u(t)}\tg)
$$ 
of the fixed manifold $f_g(M)$ with a varying metric on it. We have
then equations arising by the transformations of the exterior quantities in 
(\ref{sce}) associated to $f_{g_t}(M)$ and $f_g(M)$, respectively. 

The transformation of the exterior scalar curvature yields the equation 
\begin{equation} \label{eq1}
n(n-1) = e^{-2u(t)}(
n(n-1) 
-2(n-1)( {\rm div}_{f_g(M)} \nabla^{\tg} u^\tau- \tg(H_{f_g}, \nabla^{\tg}
u^\nu)-\|du^\tau\|^2+\frac{n}{2}\|du\|_{\tg}^2 )) \, , 
\end{equation}
where the upper indices $\tau$ and $\nu$ stand for the tangential and normal 
components, respectively, while the transformations of the squared norm of 
the mean curvature vector and second fundamental forms on  
$f_{g_t}(M)$ and $f_g(M)$ yield the equations 
\begin{eqnarray}
\| H_{f_{g_t}} \|^2 & = & e^{-2u(t)}(\| H_{f_g}\|^2 
- 2n\tg(H_{f_g},\nabla^{\tg}u ^\nu)
+n^2 \tg(\nabla^{\tg}u ^\nu, \nabla^{\tg} u^\nu)) \, , \label{eq3}\\
\| \alpha_{f_{g_t}}\|^2 & = & e^{-2u(t)}(\| \alpha_{f_g}\|^2 - 2\tg(H_{f_g},
\nabla^{\tg}u^\nu) +n \tg(\nabla^{\tg}u ^\nu, \nabla^{\tg} u^\nu)) \, ,   
\label{eq4} 
\end{eqnarray}
respectively. By (\ref{sce}), we then have that
\begin{equation} \label{eq5}
s_{g_t}= e^{-2u(t)}\left( s_{g}-2(n-1){\rm div}_{f_g(M)}(\nabla^{\tg}u^\tau)
-(n-1)(n-2)\tg(\nabla^{\tg}u^{\tau},\nabla^{\tg}u^\tau)  \right) \, ,   
\end{equation}
and that the densities of the functionals $\mc{W}_{f_{g_t}}(M)$ and  
$\mc{D}_{f_{g_t}}(M)$ in (\ref{eq12}) are 
$$
\begin{array}{l}
e^{-2u(t)}( n^2+ 2n\Delta^{g}u +\|H_{f_g}\|^2+n(n-2) \|du^\tau\|^2) \, ,\\
e^{-2u(t)}( n+ 2\Delta^{g}u +\|\alpha _{f_g}\|^2+(n-2) \|du^\tau\|^2) \, ,
\end{array}
$$
respectively, so though defined through their
isometric embeddings, $\mc{W}_{f_g}(M)$ and $\mc{D}_{f_g}(M)$
are well-defined intrinsic functionals in the space 
of metrics in the conformal class of $g$.  

\subsection{Minimal embeddings of surfaces}
In the two dimensional case, 
the functional $\mc{W}_{f_g}(M)$ in (\ref{eq12}) over the sphere $\mb{S}^3 
\hookrightarrow \mb{S}^{\tn}$ background is known as the Willmore energy of 
$f_g(M)$. We rename $M^2$ as $\Sigma$, and denote by $\mc{M}_{[g]}(\Sigma)$ and
$\mc{C}(\Sigma)$ the spaces of metrics of $\Sigma$ in the conformal class 
$[g]$ and set of conformal classes of metrics, respectively. 

By (\ref{eq1}), (\ref{eq3}) and (\ref{eq4}) we obtain the identities
$$ 
\begin{array}{rcl}
\Theta_{f_{g_t}}(\Sigma) & = & \Theta_{f_g}(\Sigma)+2{\displaystyle 
\int_{f_g(\Sigma)} (\< H_{f_g},\nabla^{\tg} u^\nu\> - \< \nabla^{\tg} u^\nu, 
\nabla^{\tg} u^\nu\>)} d\mu_g \, ,  \vspace{1mm} \\
\mc{W}_{f_{g_t}}(\Sigma) = 
2\Theta_{f_{g_t}}(\Sigma)+\Psi_{f_{g_t}}(\Sigma) & = & 2\Theta_{f_{g}}(
\Sigma) + \Psi_{f_g}(\Sigma) =
\mc{W}_{f_g}(\Sigma)  \, , \vspace{1mm} \\ 
\mc{D}_{f_{g_t}}(\Sigma) = 
\Theta_{f_{g_t}}(\Sigma)+\Pi_{f_{g_t}}(\Sigma) & = & \Theta_{f_g}(\Sigma)+
\Pi_{f_g}(\Sigma)
= \mc{W}_{f_{g}}(\Sigma) \, , 
\end{array}
$$
the first one of which implies that if 
$f_{g_t}$ is conformally related to a minimal $f_g$, then 
${\rm area}_{g_t}(\Sigma)\leq {\rm area}_g(\Sigma)$, while the latter two show
the invariance of each, $\mc{W}_{f_g}(\Sigma)$ and $\mc{D}_{f_g}(M)$, under 
conformal deformations.  The total scalar curvature functional of $(\Sigma, g)$
is the topological constant 
$$
\mc{S}_g(\Sigma)= 4\pi \chi(\Sigma)= \mc{W}_{f_g}(\Sigma)-\mc{D}_{f_g}(\Sigma)
\, ,
$$
written as the difference of functionals that are invariant 
under conformal deformations $f_{g_t}$ of $f_g$. We thus have that 
if $F: (\mb{S}^{\tn},\tg) \rightarrow (\mb{S}^{\tn},\tg)$ is a conformal 
diffeomorphism and $f_g(\Sigma)\in \mc{M}_{\mb{S}^{\tn}}(\Sigma)$, then 
$F(f_g(\Sigma)) \hookrightarrow (\mb{S}^{\tn}, \tg)$ has intrinsic metric 
$\tg\mid_{F(f_g(\Sigma))}$, and $\mc{W}_{f_g}(\Sigma)= 
\mc{W}_{F\circ f_{g}}(\Sigma)$, so we naturally define
\begin{equation} \label{cci}
\mc{W}(\Sigma,[g])=\inf_{g'\in \mc{M}_{[g]}(\Sigma)}\mc{W}_{f_{g'}}(\Sigma) =
\inf_{g'\in \mc{M}_{[g]}(\Sigma)} \int_{f_{g'}(\Sigma)} 
\left( 4 + \| H_{f_{g'}}\|^2 \right) d\mu_{f_{g'}} \, ,  
\end{equation}
and attempt to find its optimal isometric embedding realizer so that
we can compare the values of the class invariant differentiable function 
$$
\mc{C}(\Sigma) \ni [g] \rightarrow \mc{W}(\Sigma,[g])
$$
across conformal classes of metrics by comparing the values that the optimal 
realizer of each class achieves. For any surface $\Sigma^k$ of 
diffeomorphism type fixed by the topological genus $k$, there exists a 
distinguished conformal class $[g_k]$,  
and metric $g_k$ representing it of minimal isometric embedding $f_{g_k}
\in \mc{M}_{\mb{S}^{\tn}}(\Sigma)$, such that
\begin{equation} \label{ew2}
\mc{W}_{f_g}(\Sigma) = \mc{W}(\Sigma,[g]) \geq \mc{W}(\Sigma,[g_k])=
\mc{W}_{f_{g_k}}(\Sigma)\, ,
\end{equation}
with the equality achieved if, and only if, $[g]=[g_k]$ and $f_{g}$ is
conformally equivalent to $f_{g_k}$ in the sphere background 
\cite[Theorems 1 \&  9]{sim2}. Notice that if $g_t$ is a path of conformal
deformations of $g$ so $[g_t]=[g]$, the Palais' path of diffeomorphism 
$F_t: \mb{S}^{\tn} \rightarrow \mb{S}^{\tn}$ for which $F_t^{*}\tg
\mid_{f_{g_t}(\Sigma)}=g_t$ starts at
$\BOne_{\mb{S}^{\tn}}$, and set theoretically fixes $f_g(\Sigma)$ for all $t$,
with the metric $F_t^{*}\tg=e^{2u(t)}\tg$ inducing the intrinsic conformal 
metric deformation of $g$ when restricted to the submanifold $f_g(\Sigma)$, 
and $\mc{W}_{f_{g_t}} (\Sigma)=\mc{W}(\Sigma,[g])$. 
Although the values of conformally invariant 
Riemannian functionals along the $g_t$s stay the same, the nontensorial 
quantities $\nabla^{g_t}$ vary with the $t$-dependent gauge over 
$f_g(\Sigma)$. If the starting metric $g$ has minimal isometric embedding
$f_g$ and the conformal deformations $g_t$ are through area preserving 
metrics, the minimality of the embeddings $f_{g_t}$ remains,            
assertion last that comes at the expense of knowing the
area of the metric on $\Sigma$ with minimal isometric embedding $f_g$, and
the finding of an appropriately related path of background sphere 
diffeomorphism $F_t$ realizing the minimal conformal deformations of $f_g$.    

\begin{theorem}
Suppose that $f_g: (\Sigma^k,g) \rightarrow (\mb{S}^{2+q},\tg)$ is a minimal
isometric embedding of a closed surface $(\Sigma^k,g)$,  
and that $[0,1]\ni t \rightarrow g_t=e^{2u(t)}g$ is a 
path of area preserving conformal deformations of $g$ on $f_{g}(\Sigma)$. 
Then there exists a path of conformal diffeomorphism $F_t: (\mb{S}^{2+p},
F_t^* \tg) \rightarrow (\mb{S}^{2+p},\tg)$ such that
$F_t\mid_{t=0}=\BOne_{\mb{S}^{2+p}}$, $F_t(f_g(M))=f_g(M)$, $F^{*}_t \tilde{g}
\mid_{f_g(M)}=g_t$ and $f_{g_t}: (\Sigma,g_t) \rightarrow (\mb{S}^{2+p},\tg)$ 
is a path of conformal deformations of $f_g$ by minimal embeddings. 
In particular, if $s_g$ is not constant, there exists a path $f_{g_t}: 
(\Sigma,g_t) \rightarrow (\mb{S}^{2+q},\tg)$  of area $\mu_g(\Sigma)$ 
minimal conformal deformations of $f_g$ to the isometric 
embedding $f_{g'}$ of a metric $g'\in \mc{M}_{[g]}(\Sigma)$ of 
scalar curvature $4\pi \chi(\Sigma)/\mu_g(\Sigma)$.    
\end{theorem}

{\it Proof}. Along any conformal deformation $f_{g_t}$ of $f_g$, 
by (\ref{eq1}) we have that
$$
2=e^{-2u}(2+2\Delta^g u -2 \tg(\nabla^{\tg}u^{\nu},\nabla^{\tg}u^{\nu}))\, ,  
$$
so if the deformation is area preserving, by integration of this identity
with respect to the area measure $e^{2u}d\mu_g$ of $g_t$, we conclude that 
$\| \nabla^{\tg}u^{\nu}\|^2 =0$, and so $\nabla^{\tg}u^{\nu}$ is the zero 
vector along points of $f_{g_t}(\Sigma) \hookrightarrow \mb{S}^{2+q}$.

We let $u=u(t)$ be the function on $f_g(\Sigma)$ defining the conformal
deformation $g_t=e^{2u(t)}g$ of $g$.  
By the tube lemma, there exist some $\varepsilon > 0$ and a smooth extension 
$u_t$ of $u(t)$ to the whole of $\mb{S}^{2+q}$ that is 
a function of the distance in the normal directions to the submanifold 
$f_g(\Sigma)$, supported on the 
$2\varepsilon$ tubular neighborhood of $f_g(\Sigma)$, and    
constant in the normal directions in the $\varepsilon$ tubular neighborhood 
inside it, so $\nabla^{\tg}u_t$ restricts tangentially to the boundary of any 
$r< \varepsilon$ tubular neighborhood of $f_g(\Sigma)$, and 
on the zero section, $\nabla^{\tg}u_t^\nu \mid_{f_g(\Sigma)}=0$.

The flow $F_t$ of the time dependent vector field $\nabla^{\tg}u_t$ 
on $\mb{S}^{2+q}$ generates 
Palais' extension diffeomorphisms  $F_t: \mb{S}^{2+q} \rightarrow \mb{S}^{2+q}$
such that such that $F^*_t\tg = e^{2u_t}\tg$ (cf. \cite[Theorem p. 936]{past}), 
and so 
$F^*_t\tg \mid_{f_{g_t}(\Sigma)}= e^{2u(t)}\tg\mid_{{f_g}(\Sigma)}=
g_t$. Since $\nabla^{\tg}u_t^\nu\mid_{f_g(\Sigma)}=0$, by (\ref{eq3}), 
$\| H_{f_{g_t}}\|^2=0$. This proves the first part of the statement.  

Suppose now that the scalar curvature $s_g$ of $(f_g(\Sigma),g=\tg
\mid_{f_g(M)})$ is not constant. By the 
uniformization theorem if the surface is orientable, or by a 
combination of that theorem with a 2-to-1 covering of orientable cover 
if the surface is not, there exists a path $u=u(t)$ of conformal factors 
$e^{2u}$, $u\mid_{t=0}=0$, such that $[0,1]\ni t \rightarrow g_t=e^{2u}g$ is
a path of area preserving metrics on $f_g(\Sigma)$ satisfying the equation
$$
s_{g_t}= e^{-2u}(s_g +2\Delta^g u) 
$$    
and $s_{g_t\mid_{t=1}}=\frac{4\pi \chi(\Sigma)}{\mu_g(\Sigma)}$. Along this
path of metrics, we have that
$$
(4\pi \chi(\Sigma))^2=\left( \int s_{g_t} d\mu_{g_t}\right)^2
 \leq \mu_g(\Sigma) \int s_{g_t}^2 d\mu_{g_t}\, , 
$$
and the path ends at the said $g_1$ by ensuring that $g_1$ minimizes the
nonconformally invariant functional given by the squared $L^2$ norm of the
scalar curvature, which is always possible, so the extreme case of the 
Cauchy-Schwarz inequality above is achieved.     
By the preliminary portion of the proof, we may then realize this 
path as the metrics induced by the associated gauge changing metric 
$F_t^* \tg$ of $\mb{S}^{2+q}$ on the minimally embedded surface 
$f_{g_t}(\Sigma)$ that set theoretically is $f_g(\Sigma)$.      
\qed 

In the case of a Riemann surface $\Sigma^k$ of genus $k$,
the distinguished conformal class in (\ref{ew2}) is
that of the equilateral Lawson Riemann surface $\xi_{k,1}$,
and we have that  
\begin{equation} \label{weci}
\mc{W}_{f_g}(\Sigma)=\mc{W}(\Sigma,[g]) \geq 
\mc{W}_{f_{g_{\xi_{k,1}}}}(\Sigma)= \mc{W}(\xi_{k,1},[g_{\xi_{k,1}}])= 
4 \mu_{g_{\xi_{k,1}}}(\xi_{k,1}) \, , 
\end{equation}
with the equality achieved if, and only if, $f_g$ is conformally equivalent
to $f_{g_{\xi_{k,1}}}$ \cite[Theorem 1]{sim2}. We then have that
$\mc{W}(\xi_{k,1},[g_{\xi_{k,1}}])$ increases monotonically with $k$ to the
upper bound $32\pi$. In the elliptic case, this result is due to
Almgren, who proved the more general version asserting that the
only minimal immersion of a genus $0$ surface into $\mb{S}^3$ must
be that of the standard sphere \cite{alm}. The 
case of $k=1$, when the optimal lower bound is achieved by the Clifford torus 
$\xi_{1,1}$, is due to Cod\'a \& Neves \cite{cone}, who proved also a more 
general statement when solving in the affirmative the original Willmore 
conjecture. These two cases can be obtained by using the
lower and upper end, respectively, of the gap theorem of Simons
\cite[Theorem 5.3.2, Corollary 5.3.2]{si} \cite[Main Theorem]{cdck} 
\cite[Corollary 2]{law0}, though the $k=1$ case follows by a general 
argument that applies for any $k\geq 1$. But its alternative proof 
using Simons' gap theorem leads in addition to a proof of the theorem 
of Brendle \cite{bren} on
an earlier conjecture of Lawson \cite{law2} stating that the  
flat torus $(\xi_{1,1},g_{\xi_{1,1}})$ is the only Riemann surface of genus 
one that can be minimally embedded into $\mb{S}^3$ \cite[Corollary 3]{sim2}. 
For the metric $g$ of any such embedded torus can be conformally deformed to 
a flat metric $g'$ of the same area; by the gap theorem, $g'=g_{\xi_{1,1}}$
up to isometries of the background; 
and by Simons' identity on the rough Laplacian of the second fundamental 
form \cite[Theorem 5.3.1]{si}, we can then show that 
$\mu_g= \mu_{g_{\xi_{1,1}}}$, and so $g=g'=g_{\xi_{1,1}}$. 

If $\Sigma^k$ is nonorientable, if $k=1$, $\Sigma^1=\mb{P}^2(\mb{R})$ and 
$[g_1]$ is the projective space only conformal class, with 
$\mu_{g_1}(\mb{P}^2(\mb{R}))=6\pi$. If $k=2$, 
$\Sigma^k=\tilde{\tau}_{3,1}$, $[g_2]=[g_{\tilde{\tau}_{3,1}}]$ and 
$\mu_{g_2}(\tilde{\tau}_{3,1})= 6\pi E(2\sqrt{2}/3)$, as observed earlier.
The remaining cases are not so explicitly known, but the areas of the
optimal metric representatives $g_k$ are known to be bounded above by $8\pi$, 
and by the explicit known values above, 
$\mc{W}(\Sigma^k,[g_k])$ is monotonically increasing with $k$ for at least
$k<3$ (see \cite[Remark 11]{sim2}). 

\begin{theorem}
Suppose that $(\Sigma,g)=(\xi_{k/m,m},g_{\xi_{k/m,m}})$ where  
$m$ any integer divisor of $k>1$. 
Then for any smooth path $[0,1]\ni t \rightarrow g_t$ of  
area $a=\mu_{g_{\xi_{k/m,m}}}(\Sigma)$ 
conformal deformations of $g_{\xi_{k/m,m}}$ in the conformal class
$[g_{\xi_{k/m,m}}]$ to a metric $g_a$ of constant  
scalar curvature $4\pi \chi(\Sigma)/a$, 
there exists an associated path $f_{g_t}: (\Sigma,g_t)\rightarrow (\mb{S}^3,
\tg) $ of minimal isometric conformal deformations of $f_{g_{\xi_{k/m,m}}}$
to the isometric embedding $f_{g_a}$ of $g_a$. 
\end{theorem}

{\it Proof}. Since $k>1$, there are points on $\Sigma$ where the metric
$g_{\xi_{k/m,m}}$ is such that $s_{g_{\xi_{k/m,m}}}=2$ \cite[Proposition
1.5]{la2}, and so, by the Gauss-Bonnet theorem, $s_{g_{\xi_{k/m,m}}}$ is not
the constant function. By the uniformization theorem, 
the metric $g_{\xi_{k/m,m}}$ on $f_g(\Sigma)=\xi_{k/m,m}$ can be deformed  
by area preserving metrics in its conformal class to a metric 
$g_a$ of constant scalar curvature, which by the Gauss-Bonnet theorem again, 
must be of scalar curvature $4\pi \chi(\Sigma)/a$. The result now follows by 
Theorem 1, as any area preserving conformal deformation of   
$g_{\xi_{k/m,m}}$ to $g_a$ is realized by 
an associated gauge changing conformal diffeomorphism 
$F_t: \mb{S}^{3} \rightarrow \mb{S}^{3}$ that fixes $f_g(\Sigma)=\xi_{k/m,m}$
and is such that $F^*_t\tg \mid_{f_g(\Sigma)}=g_t$, with corresponding family
of minimal isometric conformal embedding deformations $f_{g_t}$ of
$f_{g_{\xi_{k/m,m}}}$ to the isometric embedding $f_{g_a}$ of $g_a$.  
\qed

\begin{corollary}
If $f_{g_a}$ is the minimal isometric embedding of the equal area  
constant scalar curvature metric $g_a$ in the conformal class of 
$g_{\xi_{k/m,m}}$, then
$$
\| \alpha_{f_{g_a}}\|^2= 2-\frac{4\pi\chi(\xi_{k/m,m})}{\mu_{g_{\xi_{k/m,m}}}(
\xi_{k/m,m})}=2-\frac{8\pi(1-k)}{\mu_{g_{\xi_{k/m,m}}}(\xi_{k/m,m})} \, .
$$
\end{corollary}

{\it Proof}. This follows by (\ref{sce}) and the Gauss-Bonnet theorem. 
\qed

\begin{remark}
Besides the $\xi_{m',m}$s of Lawson, there are other families of Riemann 
surfaces embedded minimally into $\mb{S}^3$, which interestingly enough, are 
all of relatively large genus \cite{kps, chso, kaya}. With each of them 
playing the role of the $\xi_{k/m,m}$ surfaces in the statement, Theorem 2 
applies. Theorem 1 applies to any of the Willmore surfaces of Bryant 
\cite{br0} immersed into $\mb{S}^4$ when the immersion is an embedding, or 
equally well to the bipolar surface $\tilde{\tau}_{3,1}$ of topological 
genus $2$ minimally embedded into $\mb{S}^4$ with area 
$6\pi E(2\sqrt{2}/3)$, and for which $g_{\tilde{\tau}_{3,1}}$ is not itself 
of constant scalar curvature.    
\end{remark}

\section{The notion of $\sigma$ invariant for surfaces}
We let $\mc{M}_{a,[g]}(\Sigma)$ be the submanifold of metrics of area $a$ in 
the conformal class of $g$, 
set $a=a([g]):=\frac{1}{4}\mc{W}(\Sigma,[g])$, 
and consider the functional 
$$
\mc{S}^2_a: \mc{M}_{a,[g]}(\Sigma)\ni g \rightarrow \int s_g^2 d\mu_g \, . 
$$
This is not a conformally invariant functional, 
but its minimum is always achieved by a metric of constant scalar curvature
$4\pi \chi(\Sigma)/a$,
 and for dimensional reasons, this minimum is the end point of a path of area 
preserving conformal deformations of $g$.
By the Gauss-Bonnet theorem, and the extreme case of the 
Cauchy-Schwarz inequality, the critical value at this minimizer is 
$$
\mc{S}^2_a(\Sigma,[g]) = \frac{(4\pi \chi(\Sigma))^2}{a} \, . 
$$
The canonical realizer of the infimum (\ref{cci}) defining 
$\mc{W}(\Sigma,[g])$ is thus the $f_g$ of an Einstein metric $g$ in $[g]$ of 
minimal area \cite[Theorem 6]{sim2}, 
 which if $\chi(\Sigma) \leq 0$, can be uniquely associated
to the fundamental domain of the conformal class $[g]$ when 
$\Sigma$ is orientable, or to that of the conformal class of the lifted metric
to its 2-to-1 orientable cover when $\Sigma$ is not \cite[Theorem 8]{sim2}. 
The distinguished conformal class $[g_k]\in \mc{C}(\Sigma)$ in 
(\ref{ew2}) is the one whose associated fundamental domain is the most 
symmetric of them all, and the associated metric representative $g_k$ has 
the smallest area among the canonical realizers of the various elements 
of $\mc{C}(\Sigma)$, fixing the optimal scale among them.

Thus, we define  
$$
\sigma^2(\Sigma) = \sup_{[g]\in \mc{C}(\Sigma))} 
\mc{S}^2_{a([g])}(\Sigma,[g])\, .   
$$
In the orientable case, the bounds (\ref{weci}) yield the following:

\begin{theorem}
If $\Sigma=\Sigma^k$ is a Riemann surface of genus $k$, and
$f_g: (\Sigma,g) \rightarrow (\mb{S}^{\tn},\tg)$, we have that
$$
\frac{(4\pi \chi(\Sigma))^2}{\frac{1}{4}\mc{W}_{f_g}(\Sigma)} 
\leq \sigma^2(\Sigma)= 
\frac{(4\pi \chi(\Sigma))^2}{\frac{1}{4}\mc{W}(\Sigma, [g_{\xi_{k,1}}])}
=\frac{(4\pi \chi(\Sigma))^2}{\mu_{g_{\xi_{k,1}}}(\xi_{k,1})}
 \, , 
$$   
and for $k\neq 1$, the equality is achieved by $f_g$ if, and only if, 
$[g]=[g_{\xi_{k,1}}]$ and $f_g$ is conformally equivalent to 
$f_{g_{\xi_{k,1}}}$, while for $k=1$, if $\mc{W}(\Sigma,[g])=\mc{W}(
\Sigma,[g_{1,1}])$ then $[g]=[g_{1,1}]$ and $f_g$ is conformally
equivalent to $f_{g_{1,1}}$.     
\end{theorem}

In the general situation, irrespective of the orientation, of 
a surface $\Sigma^k$ of topological genus $k$ and distinguished conformal 
class $[g_k]$ for which (\ref{ew2}) holds, we have that
$$
\frac{(4\pi \chi(\Sigma))^2}{\frac{1}{4}\mc{W}_{f_g}(\Sigma)}
\leq \sigma^2(\Sigma)= \frac{(4\pi \chi(\Sigma))^2}{\frac{1}{4}\mc{W}(
\Sigma, [g_k])} \, ,
$$
and if $\chi(\Sigma) \neq 0$, the equality is achieved by $f_g$ if, and only if,
$[g]=[g_k]$ and $f_g$ is conformally equivalent to $f_{g_k}$, while if
$\chi(\Sigma)=0$ and $\mc{W}_{f_g}(\Sigma)=\mc{W}(\Sigma,[g_k])$ 
then $[g]=[g_{\xi_{1,1}}]$ if $\Sigma^k$ is orientable or 
$[g]=[g_{\tilde{\tau}_{3,1}}]$ otherwise, and in either case, 
$f_g$ is conformally equivalent to $f_{g_k}$.  
Thus, we always have an optimal conformal class that
achieves $\sigma^2(\Sigma)$, even when $\chi(\Sigma)=0$. We set 
$$ 
\sigma(\Sigma,[g])=\frac{8\pi \chi(\Sigma)}{
\mc{W}(\Sigma,[g])^{\frac{1}{2}}}\, ,
$$ 
and define the sigma invariant of the surface by
$$ 
\sigma(\Sigma) =
\sup_{[g]\in \mc{C}(\Sigma)} 
\sigma(\Sigma,[g])
={\rm sign}(\chi(\Sigma))\sqrt{\sigma^2(\Sigma)}=
\frac{8\pi \chi(\Sigma)}{\mc{W}(\Sigma,[g_k])^{\frac{1}{2}}}\, ,
$$ 
taking $f_{g_k}$ as its optimal realizer.  
The sign of $\chi(\Sigma)$ determines the sign of any constant scalar curvature
metric $g$ in $\mc{M}(\Sigma)$, regardless of scale or conformal class, and 
we have the universal bound 
$$
\sigma(\Sigma) \leq \sigma(\mb{S}^2) = 4\pi^{\frac{1}{2}} \, .
$$
 
\begin{theorem}
Let $\Sigma^k$ be a surface of topological genus $k$ and distinguished 
conformal class $[g_k]$. Then $16\pi \leq \mc{W}(\Sigma,[g_k]) < 32\pi$, 
and   
$$ 
\mc{W}(\Sigma,[g_k])\leq \mc{W}(\Sigma,[g]) 
$$
for any conformal class of metrics $[g]$ on $\Sigma$.  
\begin{enumerate}[label={\rm \alph*)}] 
\item If $\chi(\Sigma) < 0$, then 
$\sigma(\Sigma,[g_k]) \leq \sigma(\Sigma,[g])$.
\item If $\chi(\Sigma) \geq 0$, then $\sigma(\Sigma,[g]) \leq 
\sigma(\Sigma,[g_k])$. 
If $\chi(\Sigma)=0$, the distinguished conformal classes in the orientable and
nonorientable cases are $[g_{\xi_{1,1}}]$
and $[g_{\tilde{\eta}_{3,1}}]$ with Willmore energies  
$8\pi^2$ and $24\pi E(2\sqrt{2}/3)$, respectively, while 
if $\chi(\Sigma)>0$, $\Sigma$ orientable or not carries a unique conformal 
class of metrics, $\sigma(\mb{S}^2)= 4 \pi^{\frac{1}{2}}$ and the Willmore 
energy of its class is $16\pi$, while 
$\sigma(\mb{P}^2(\mb{R}))=\frac{4}{\sqrt{6}}\pi^{
\frac{1}{2}}$ and the Willmore energy of its class is $24\pi$. 
\end{enumerate}
\end{theorem}


\begin{thebibliography}{69}
\bibitem{alm}
F.J. Almgren, {\it Some interior regularity theorems for minimal surfaces
and an extension of Berstein's theorem}. Ann. of Math. 84 (1966), pp. 277-292.
\bibitem{bren}
S. Brendle, {\it Embedded minimal tori in $\mb{S}^3$ and the Lawson 
conjecture}, Acta Math., 211 (2013), pp. 177-190.
\bibitem{br0}
R. Bryant, {\it Conformal and minimal immersions of compact surfaces into the
4-sphere}, J. Diff. Geom. 17 (1982), pp. 455-473.
\bibitem{cdck}
S.S. Chern, M. do Carmo \& S. Kobayashi, {\it Minimal submanifolds of a 
sphere with second fundamental form of constant length}.  1970  
Functional Analysis and Related Fields (Proc. Conf. for M. Stone, Univ. 
Chicago, Chicago, Ill., 1968)  pp. 59-75 Springer, New York.
\bibitem{chso}
J. Choe \& M. Soret, {\it New minimal surfaces in $S^3$ desingularizing the   
Clifford tori}, arXiv:1304.3184
\bibitem{cone}
F. Cod\'a Marques \& A. Neves, {\it Min-max theory and the Willmore 
conjecture}, Ann. of Math., 179 (2014), pp. 683-782. 
\bibitem{kaya}
N. Kapouleas \& S.D. Yang, {\it Minimal surfaces in the three-sphere by 
doubling the Clifford torus}, Amer. J. Math. 132 (2010), pp. 257–295.
\bibitem{kps} 
H. Karcher, U. Pinkall \& I. Sterling, {\it New minimal surfaces in $S^3$}, 
J. Diff. Geom. 28 (1988), pp. 169–185.
\bibitem{ku}
R. Kusner, {\it Comparison surfaces for the Willmore problem}, Pac. J. Math. (2)
138 (1989), pp. 317-345.
\bibitem{kls}
E. Kuwert, Y. Li \& R. Sch\"{a}tzle, {\it The large genus limit of the infimum 
of the Willmore energy}, Amer. J. Math. 132 (2010), pp. 37–51.
\bibitem{lap}
H. Lapointe, {\it Spectral properties of bipolar minimal surfaces in 
$\mb{S}^4$}, Diff. Geom. Appl. 26 (2008), pp. 9-22. 
\bibitem{la2}
H.B. Lawson Jr. {\it Complete minimal surfaces in $\mb{S}^3$}, 
Ann. of Math. (2) 92 (1970), pp. 335-374.
\bibitem{law0}
H.B. Lawson Jr. {\it Local rigidity theorems for minimal hypersurfaces},
Ann. of Math. (2) 89 (1969), pp. 187-197.
\bibitem{law2}
H.B. Lawson Jr. {\it The unknottedness of minimal embeddings}, 
Inv. Math. (11) (1970), pp. 183-187.
\bibitem{nash}
J. Nash, {\it The imbedding problem for Riemannian manifolds}, Ann. of Math.
63 (1956), pp. 20-63.
\bibitem{pal}
R.S. Palais, {\it Local triviality of the restriction map for embeddings},
Comment. Math. Helv. 34 (1960), pp. 305–312.
\bibitem{past}
R.S. Palais \& T.E. Stewart, {\it Deformations of compact differentiable 
transformation groups}, Amer. J. Math. 82 (1960), pp. 935-937.
\bibitem{si3}
S.R. Simanca, {\it Canonical isometric embeddings of projective spaces into 
spheres}, Bol. Soc. Mat. Mexicana, 26 (2020), pp. 757-763.
\bibitem{sim2}
S.R. Simanca, {\it Conformal invariants of isometric embeddings of the 
smooth metrics of a surface}, Bull. Pol. Acad. Sci. Math. 72 (2025), pp. 
161-187.
\bibitem{scs}
S.R. Simanca, {\it $\mb{S}^6$ {\rm (}or any of $\mb{S}^2 \times \mb{S}^4$, 
$\mb{S}^2\times\mb{S}^6$, or $\mb{S}^6\times \mb{S}^6$, respectively{\rm )} is 
not  diffeomorphic to a complex manifold}, preprint 2022, arXiv:2204.12628.   
\bibitem{si}
J. Simons, {\it Minimal varieties in Riemannian manifolds}, Ann. Math. 2 
(1968), pp. 62-105. 
\end{thebibliography}
\end{document}